\documentclass[10pt,notitlepage]{article}
\usepackage{epsfig}
\usepackage{amssymb}

\newcommand{\be}{\begin{equation}}
\newcommand{\ee}{\end{equation}}
\newcommand{\bea}{\begin{eqnarray}}
\newcommand{\eea}{\end{eqnarray}}

\usepackage{bm}
\usepackage{color}
\usepackage{amsmath}
\usepackage{graphicx}
\usepackage{amssymb}
\usepackage{cite}

\begin{document}


\title{A fast memoryless predictive algorithm in a chain
of recurrent neural networks}
\author{Boris  Rubinstein,
\\Stowers Institute for Medical Research
\\1000 50$^{}\mbox{th}$ St., Kansas City, MO 64110, U.S.A.}
\date{\today}
\maketitle

\begin{abstract}
In the recent publication \cite{Rub2020a} a fast prediction
algorithm for a single recurrent network (RN) was suggested. In this manuscript we 
generalize this approach to a chain of RNs and show
that it can be implemented in natural neural systems.
When the network is used recursively to predict
sequence of values the proposed algorithm does not require 
to store the original input sequence. It increases robustness of the new approach 
compared to the standard moving/expanding window predictive procedure.
We consider requirements on trained networks that allow to
implement the proposed algorithm and discuss them in 
the neuroscience context.
\end{abstract}


\section{Introduction}

Recurrent networks (RNs)  due to
their ability to process sequences of data are used in many fields 
of science, engineering and humanities, including speech, handwriting and 
human action recognition, automatic translation, robot control 
as a tools of time series prediction,
text and image generation and more complex problems in neurolinguistic
programming.
There are several basic RN architectures of different complexity 
and multiple variants of these types were discussed recently.
Multitude of recurrent network versions allows flexible combinations of different RNs in
a single complex network designed for a specific task.

After network is trained one obtains a predictive tool that should be 
applied properly to perform the desired task, i.e., one has to have
a reliable predictive algorithm to be used with the trained network.
The standard "moving window" (MW) algorithm has an input of an ordered sequence of the 
elements of similar structure, 
transforms each input element into a network state and then send the last
state into a predictor that generates a new element of the same structure as input
ones. Then a new input sequence is formed by addition of this predicted element to original input sequence 
(with its first term dropped), so that the input length remains constant and the prediction round repeats 
several times to generate more elements. 

The author of this manuscript recently suggested in \cite{Rub2020a}  
a novel  fast predictive algorithm for a system made of a single
recurrent network and a predictor. This approach does not 
require to store a part of the initial input in a short term memory
for the second and subsequent
prediction rounds. Instead the network after the first prediction
round uses its own dynamics to 
extrapolate the input. It was shown that in the neuroscience
context this memoryless (ML) algorithm is more robust and 
provides a significant speed up compared to MW approach.
The same time for an input of large enough length the extrapolations 
by both algorithms coincide.

The number of neuron connections (parameters) of RN scales as
a square of neuron number. To increase network robustness it
might be useful to replace a single RN with large number of neurons by
a chain of a few smaller RNs. 
In this manuscript we generalize ML algorithm for 
a case of a neural network made of a RN chain and a predictor.
We show that such general ML algorithm is indeed faster than 
the standard one and it can be implemented in 
natural neural systems.
It appears also that ML algorithm can be successfully applied for the 
well trained network (i.e., the network for which the deviation of the
predicted value from the ground truth value is negligibly small) that 
has high importance for neuroscience.

\section{Signal transformation in a chain of recurrent networks}
\label{algorithm0}

\begin{figure}[h!]
\begin{center}
\begin{tabular}{lc}
{\bf a} &
\psfig{figure=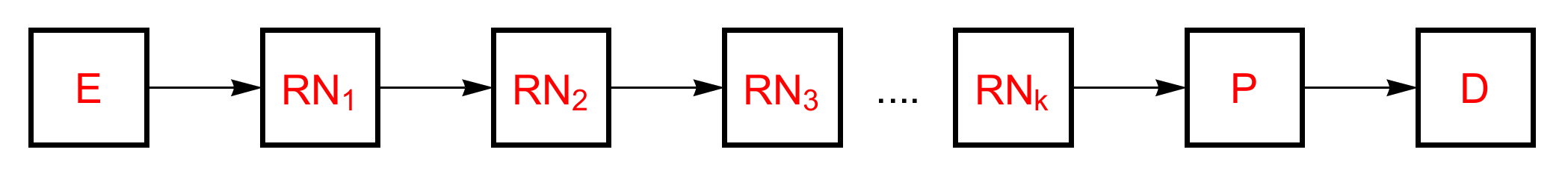,width=14.0cm} \\
{\bf b} &
\psfig{figure=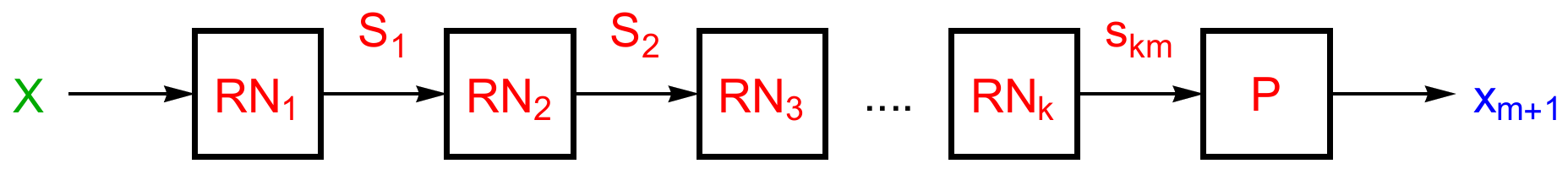,width=14.0cm} 
\end{tabular}
\caption{
({\bf a}) The architecture of a  general predictive network consisting of 
an encoder (E), a chain of  recurrent networks (RN${}_r$),
 a predictor (P) and a decoder (D). The prediction itself 
is performed in the reduced network shown in ({\bf b}). 
The input sequence $\bm X$ of the length $m$  is fed into
a chain of RNs leading to generation of the corresponding 
states sequences $\bm S_r$ (each one being an input sequence for the 
subsequent RN). The final state of the last $k$-th RN $\bm s_{k,m}$ is fed into 
the predictor P which generates a predicted vector $\bar{\bm x}_{m+1}$ of the length 
$n_0$ equal to the dimension of input vectors in  $\bm X$.
}
\label{Fig0}
\end{center}
\end{figure}

Consider a general predictive network consisting of 
an encoder, a chain of $k$ recurrent networks (RNs), a predictor and a decoder
(Fig. \ref{Fig0}a). The signal is represented by a finite length ($m$) sequence $\bm E$ of 
objects $\bm e_i, \ 1\le i \le m$ of similar structure (words, symbols, images, musical notes etc.)

Each input object $\bm e_i$ is transformed by the encoder
into a vector $\bm x_i = \bm s_{0,i}$ of the length $n_0$
thus forming an input sequence $\bm X$.
The decoder transforms a $n_0$-dimensional vector
back into an object $\bm e_{m+1}$ having 
the same structure as the input elements $\bm e_i$, i.e., 
predicts a new element based on the  sequence $\bm E$.

The heart of the predictive network is a chain between the encoder and the decoder
(Fig. \ref{Fig0}b). The input sequence $\bm X$ of the length $m$  is fed into
a chain of RNs leading to generation of the corresponding 
states sequences $\bm S_r$ (each one being an input sequence for the 
subsequent RN). The elements of $\bm S_r$ are $n_r$ dimensional
vectors $\bm s_{r,i}$ representing an inner state
of $r$-th RN (made of $n_r$ neurons). 
Before the first element $\bm x_1$ is fed into the network the initial state
$\bm s_{r,0}$ of each RN is assumed to be a zero vector of the corresponding length
$n_r$.
The transformation of the state vector of the $r$-th RN is given by 
\be
\bm s_{r,i} = \bm F_r(\bm s_{r-1,i}, \bm s_{r,i-1}),
\label{F_r}
\ee
which describes a simple rule -- the current inner state of the RN depends on the 
previous inner state and the current input signal. This rule corresponds to an
assumption that the neural network does not store its state but just updates it
with respect to the submitted input signal and its previous state.
The final state $\bm s_{k,m}$ of the last $k$-th RN  is fed into 
the predictor P that generates a vector 
\be
\bar{\bm x}_{m+1}=\bm p(\bm s_{k,m})
\label{predictor}
\ee
having the length 
$n_0$ equal to the dimension of input vectors in  $\bm X$.
Further we will consider only the architecture shown in Fig. \ref{Fig0}b as it completely determines
the predictive properties of the network shown in Fig. \ref{Fig0}a.

For the network training one feeds it with an
input sequences $\bm X^k$ to generate a predictions $\bar{\bm x}_{m+1}^k$ which 
compared to the actual values $\bm x_{m+1}^k$. The networks parameters
are fitted to minimize the mean square difference after $N$ training rounds
$$
(1/N)\sum_{k=1}^N (\bar{\bm x}_{m+1}^k-{\bm x}_{m+1}^k)^2.
$$

\section{Predictive algorithms for trained network}
Once the network is trained it can be used to generate several
consecutive values  ${\bm x}_{m+j}, \ 1 \le j \le p$. Here and further
we drop the bar over the predicted values notation. 

\subsection{Moving/expanding window algorithm}
\begin{figure}[h!]
\begin{center}
\psfig{figure=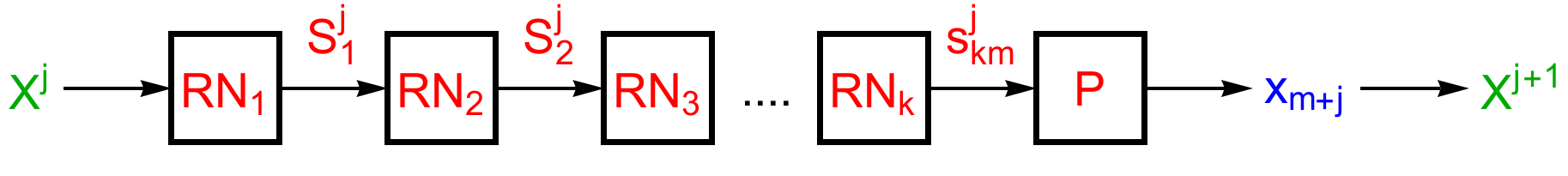,width=14.0cm} 
\caption{
The moving window algorithm prediction procedure.
The input sequence $\bm X^j$ of the length $m$  is fed into
a chain of RNs leading to generation of the corresponding 
states sequences $\bm S_r^j$ (each one being an input sequence for the 
subsequent RN). The final state of the last $k$-th RN $\bm s_{km}^j$ is fed into 
the predictor P which generates a predicted vector ${\bm x}_{m+j}$  used to 
update the input sequence $\bm X^{j+1}$ for the next prediction round. 
}
\label{MW_regular}
\end{center}
\end{figure}

When the network is trained on the sequences of the fixed length the standard
predictive algorithm uses a "moving window" (MW) recursion. 
One starts with a sequence $\bm X^{1}$ of length $m$ supplied as an input
to the network; it leads to generation of the state sequences 
$\bm S_r^j$ and the last element $\bm s_{km}^1$ of the last
state array is sent to the predictor to  
compute a prediction of the next  point  
${\bm x}_{m+1} = \bm p(\bm s_{km}^1)$. The next input sequence $\bm X^{2}$ 
is produced by dropping the first point of $\bm X^{1}$ and 
adding the predicted point ${\bm x}_{m+1}$ to the result.
This sequence is used as a new input leading to generation of ${\bm x}_{m+2}$
and a next input ${\bm X}^{3}$ is formed.
Thus at $j$-th predictive step the input $\bm X^j$ to RNN is 
formed by adding to the original input $\bm X^1$ all previously $j-1$
predicted points and shifting a "window" by $j-1$ steps forward 
(Fig. \ref{MW_regular}).
In the end of each predictive round the initial state vectors $\bm s_{r,0}^j$
of each RN either retain their values from the previous round 
(i.e., $\bm s_{r,0}^j$ = $\bm s_{r,m}^{j-1}$) or set to zero 
($\bm s_{r,0}^j$ = $\bm 0$). 
The result of prediction is not affected significantly by either choice.

The recursive procedure is repeated $p$ times to produce a sequence of  $p$ points
${\bm x}_{m+j}, (1 \le j \le p)$ approximating the sequence
$\{{\bm x}_i\}$ for $m+1 \le i \le m+p$. 
The total  number $N_{MW}^k$ of transformations (\ref{F_r}) and (\ref{predictor}) required 
to produce $p$ points is equal to $N_{MW}^k = (mk+1)p$.

For network training with sequences of the variable length 
(not larger than $M$) the 
MW predictive algorithm can be modified into the "expanding window"
(EW) version. After each predictive round the newly generated
point is added to the input sequence, so that after the $j$-th prediction 
round the length of the input sequence $\bm X^j$ 
is $m+j$, where $m$ denotes the size of the 
initial input sequence and $m+j \le M$. The main reason of EW algorithm application is that
a gradual increase of the input length usually
leads to better prediction quality.
The total  number $N_{EW}^k$ of transformations required 
to predict $p$ points is equal to 
$N_{EW}^k = N_{MW}^k + kp(p-1)/2$.
We observe that EW method requires more memory to store the input values
while MW needs a fixed memory size, so that further we focus on MW approach as more 
economical one.

In the neuroscience context an implementation of both "window" prediction algorithms
in a natural brain environment requires satisfaction of several conditions.
First, either all (for the EW) or a part (for MW) elements of the 
initial input $\bm X^1$ should be stored and reused for the second and subsequent prediction
rounds. This means that some sort of short term memory should be employed.
The neuron activation corresponding to these values should be maintained constant 
for the duration of $p$ predictive rounds. 
Second, the appearance of these values in the input sequence $\bm X^j$ should 
follow the original order of $\bm X^1$. The author of this manuscript showed \cite{Rub2020a}  that
when the first condition is not met the quality of prediction goes down 
while if the second requirement fails the accurate prediction becomes impossible.
These considerations encouraged the author to look for
an alternative predictive algorithm that does not require the input
sequence storage and thus increases predictive robustness.

\subsection{Memoryless algorithm for single RN network}
\label{single}

\begin{figure}[h!]
\begin{center}
\begin{tabular}{cc}
\psfig{figure=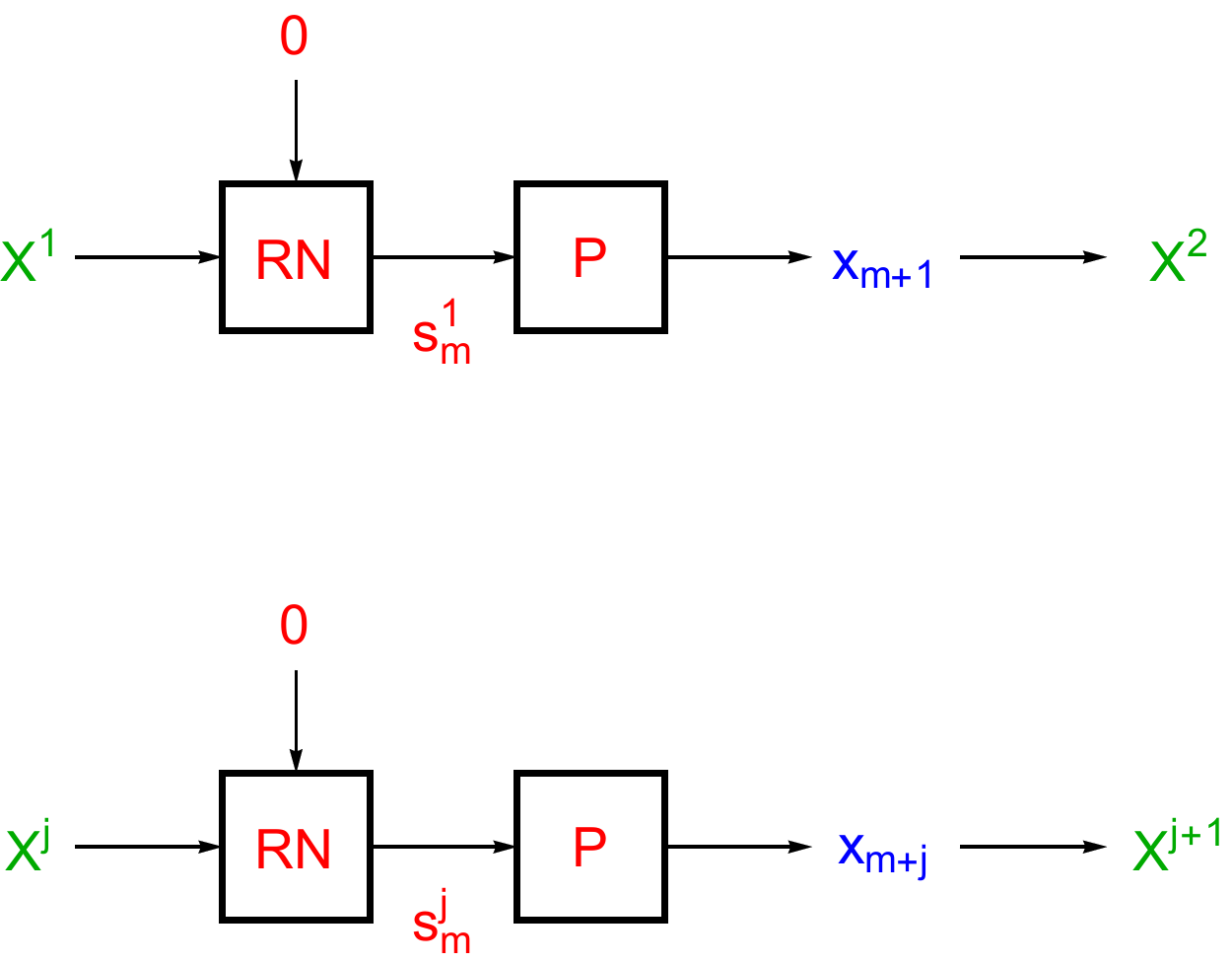,width=7.0cm}
& 
\psfig{figure=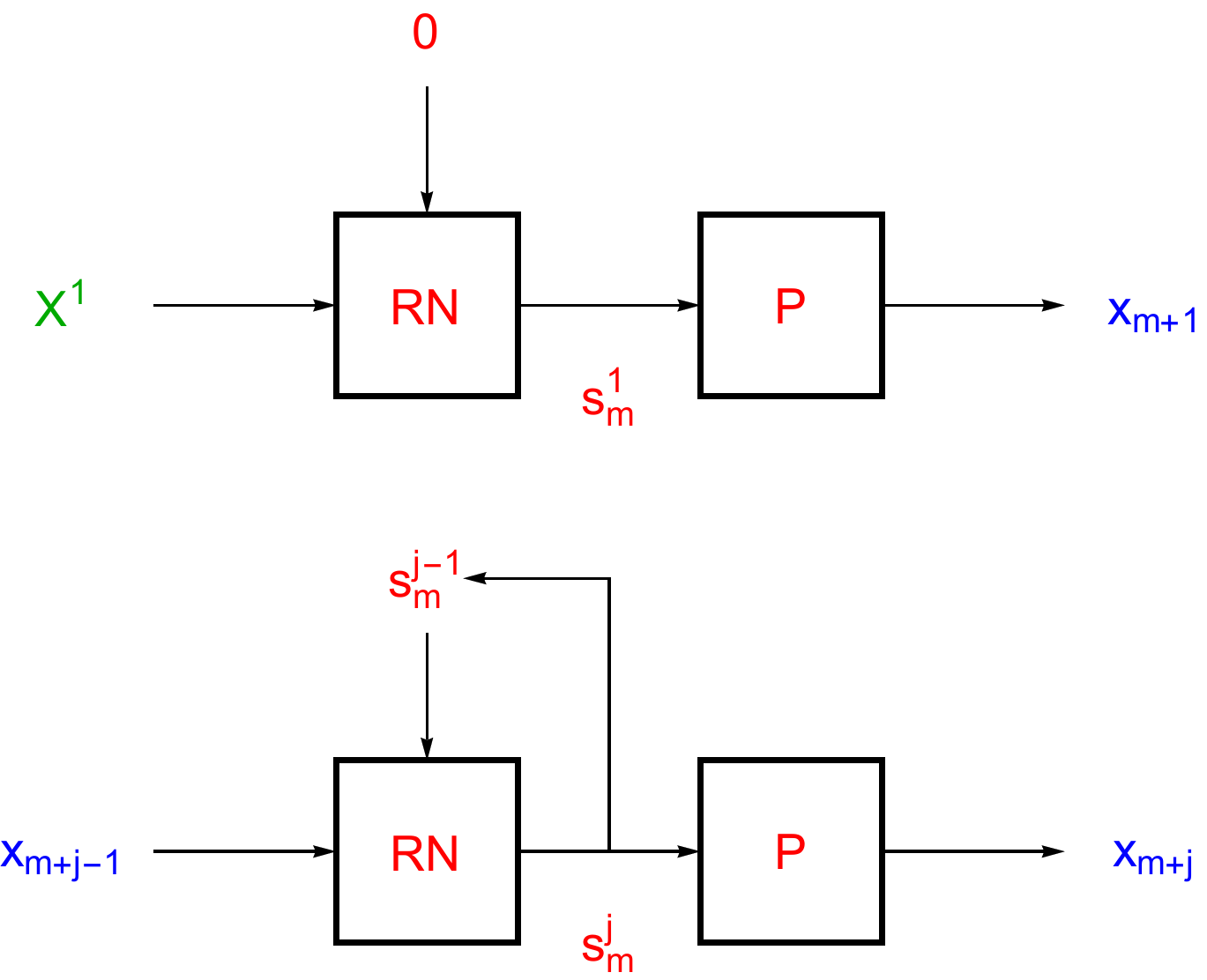,width=7.2cm} \\
{\bf a}  
& {\bf b} 
\end{tabular}
\caption{The predictive algorithms in the network consisting of 
a single recurrent network RN and
 a predictor P.
({\bf a}) The MW algorithm first (top panel) and 
a subsequent $j$-th ($j > 1$) predictive round (bottom panel). 
The input sequence $\bm X^1$ of the length $m$  is fed into
the RN with the zero initial state producing the 
state sequence $\bm S^j$ and its final state $\bm s_{m}^1$ is fed into 
the predictor P which generates a predicted vector ${\bm x}_{m+1}$ used 
for the update of the input sequence $\bm X^2$. 
This procedure is repeated at each subsequent prediction round.
({\bf b}) The memoryless (ML) algorithm first  (top panel) and 
a subsequent $j$-th ($j > 1$) predictive round (bottom panel).
The first round is similar to the one in the MW algorithm (top panel in ({\bf a}))
All subsequent rounds use the transformation (\ref{F}) of the single input value
$\bm x_{m+j-1}$ (it is the prediction made in the previous round) 
with the initial RN state $\bm s_m^{j-1}$ (it is the final RN state in the previous round).
The new value $\bm s_m^{j}$ is used as the RN state in the next round.
}
\label{singleRN}
\end{center}
\end{figure}

Consider the simplest predictive network made of a single RN ($k=1$) and a predictor.
In this case the transformations (\ref{F_r}) reduces to
\be
\bm s_{i} = \bm F(\bm x_{i}, \bm s_{i-1}),
\label{F}
\ee
where $\bm x_i$ denotes an element of the input sequence.
The MW prediction scheme is shown in Fig. \ref{singleRN}a; the top panel
shows the first predictive round while the bottom panel corresponds to 
one of the subsequent rounds.

\begin{figure}[h!]
\begin{center}
\psfig{figure=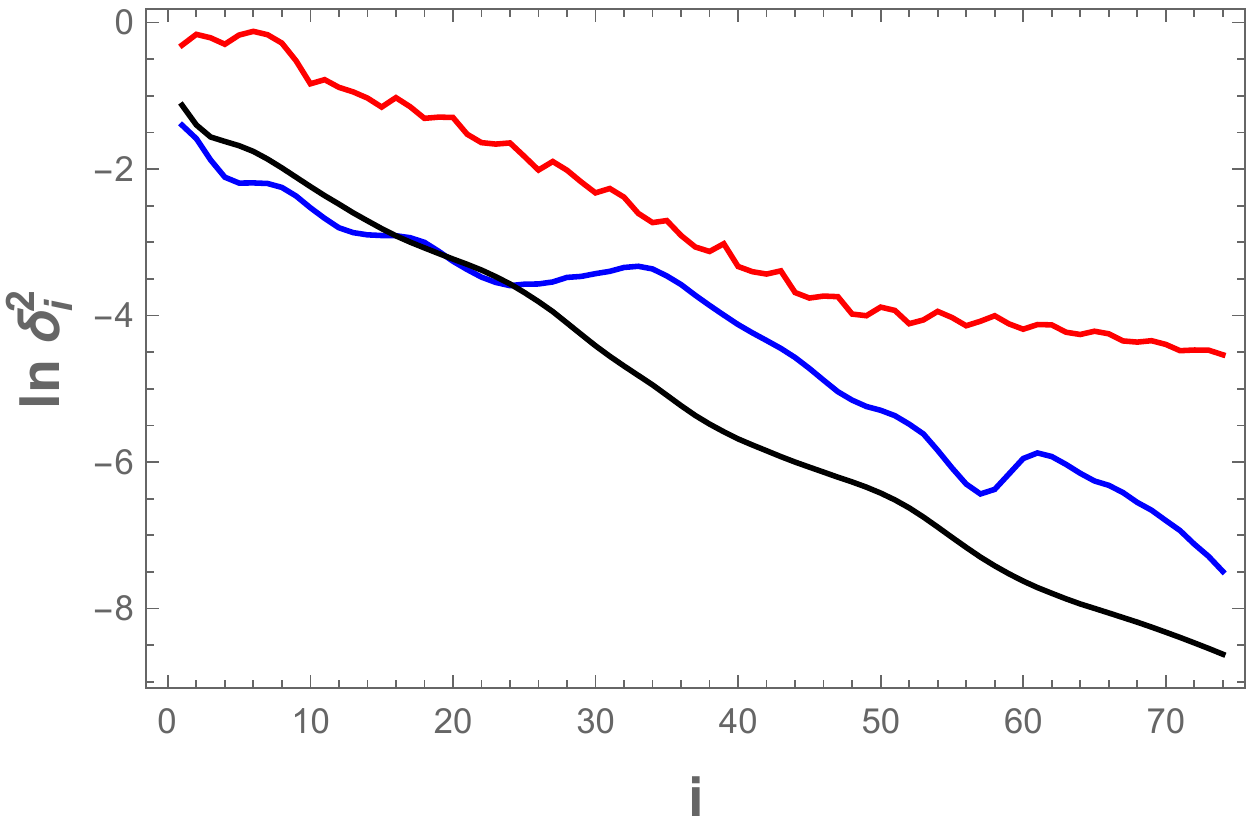,width=8.0cm}
\caption{The dynamics $\delta_{i}^2$
for the individual basic (red), gated (black) and LSTM (blue) recurrent networks with
$n=20$ neurons.}
\label{delta_dyn}
\end{center}
\end{figure}

Consider the dynamics of the RN state vectors $\bm s_{i}^j$ and  $\bm s_{i}^{j+1}$
in two adjacent prediction rounds.
First note that the MW algorithm implies the following simple relation
$\bm x_i^j = \bm x_{i-1}^{j+1}$
between the elements of input sequences $\bm X^j$ and $\bm X^{j+1}$.
It is reasonable to compare the states $\bm s_{i+1}^j$ and  $\bm s_{i}^{j+1}$
and introduce a {\it shifted difference}
\be
\bm \delta_{i}^j = \bm s_{i+1}^j-\bm s_{i}^{j+1}.
\label{delta}
\ee
It was shown in \cite{Rub2020a}  that the norm $\delta_{i}^j$
of the shifted difference decays exponentially with the transformation 
step $i$ (Fig. \ref{delta_dyn}) and for $i=m-1 \gg 1$ we find 
$\delta_{m-1}^j = |\bm s_{m}^j-\bm s_{m-1}^{j+1}| \ll 1$.
This means that the state vectors 
$\bm s_{m}^j$ and  $\bm s_{m-1}^{j+1}$ are nearly equal so that
\be
{\bm s}^{j+1}_{m-1} = {\bm s}^{j}_{m} + \bm\epsilon^j,
\quad
\epsilon^j \ll 1.
\label{last_state_diff}
\ee
This relation together with (\ref{predictor},\ref{F}) 
leads to 
\be
\bm s_m^{j+1} = \bm F(\bm x_m^{j+1},\bm s_{m-1}^{j+1})=
\bm F(\bm x_{m+j},\bm s_{m}^{j})  = 
\bm F(\bm p(\bm s_{m}^{j}),\bm s_{m}^{j}).
\label{ML_transform}
\ee
This transformation allows to compute recursively the sequence of the 
final RN states $\bm s_{m}^{j}$ for $2 \le j \le p$ by feeding the predicted value
$\bm x_{m+j-1}$ from the previous round into the 
RN with the initial state $\bm s_{m}^{j-1}$ inherited
from the same previous round.
Parallel to this computation of RN states one also fnds all the predicted values
$\bm x_{m+j}$ for $1 \le j \le p$.

The above result paves way to a memoryless (ML) predictive algorithm
replacing both MW and EW approaches.
The initial predictive round is the same as in MW/EW algorithm 
(the top panel in Fig. \ref{singleRN}b). The remaining prediction rounds
use the relation (\ref{ML_transform}) illustrated in the bottom panel in 
Fig. \ref{singleRN}b. An important feature of this algorithm that it uses
the initial sequence $\bm X^1$ only once and further it relies on
the trained network own dynamics
determined by the last expression in (\ref{ML_transform}).
The total number of transformations is 
$N_{ML}^1 = m+2p-1$ compared to $N_{MW}^1 = (m+1)p$ and
thus this approach is faster and more reliable than the standard 
MW prediction. 
The speed gain reads 
\be
\gamma_1 = \frac{N_{MW}^1}{N_{ML}^1} = 
\frac{(m+1)p}{m+2p-1}.
\label{gamma1}
\ee
Moreover, its implementation in a natural neuronal network is
simple -- as soon as the first prediction value $\bm x_{m+1}$ is obtained
it is immediately fed back into the network which internal state $\bm s_m^1$
is inherited from the first round. The result of transformation  (\ref{ML_transform})
is an updated  internal state $\bm s_m^2$ that allows to predict $\bm x_{m+2}$,
and then this process repeats.

The comparison of the predictive quality of two approaches performed in 
 \cite{Rub2020a} showed that for relatively small $m \sim 20-40$ the
ML method prediction slightly deviates from the MW extrapolation,
but for $m > 60-70$ both predictions coincide.
As the number $p$  of predicted points is usually less or approximately equal to the length $m$
of the initial input sequence the speed gain estimates is in the range from
$\gamma_1 \sim p$ for $p \ll m$ to $\gamma_1 \sim m/3$ for $p \sim m$.

\section{ML algorithm for RN chain}
\label{multi}
The numerical experiments in \cite{Rub2020a} show that for simple 
trajectories like noisy sine or triangle wave a RN with $n=10-20$ neurons
using the ML algorithm demonstrates quite high predictive quality.
For more complex cases one needs a RN with much larger number $n$ of neurons
but the number of trainable parameters scales as $n^2$ so that 
the RN robustness decreases. From this perspective it might be useful to
use instead a chain of $k > 1$ RN with $n_r$ neurons in  $r$-th RN
where 
$\sum_{r=1}^k  n_r = n$.

It is instructive to design a generalization of the ML approach to
the case of RN chain. The necessary condition of such generalization 
is a validity of an assumption that for each $r$-th RN in the chain
the shifted difference norm $\delta_{ri}^j$ decays exponentially with the 
transformation step $i$. 
To test it we consider a chain of $k=3$ RNs of different types --
the first RN is the basic recurrent network with $n_1=10$, the 
next is the LSTM network \cite{Hochreiter1997} with $n_2=15$ and the last one is
the gated recurrent network \cite{Chung2014} with $n_3=8$.

The chain is trained on the set of one-dimensional ($n_0=1$) values 
representing a sine wave $\sin(2\pi t) + a\xi(t)$ with added white noise component
with amplitude $a=0.15$.
The time step $\Delta t$ between the adjacent time points is selected equal to
$\Delta t = 0.01$. The training set consists of 12000 segments of variable
length $m$ in the range $5 \le m \le 150$.
The RNNs are trained using the Adam algorithm for 50 epochs on $80\%$ of the complete 
training set with $20\%$ validation set.
\begin{figure}[h!]
\begin{center}
\begin{tabular}{cc}
\psfig{figure=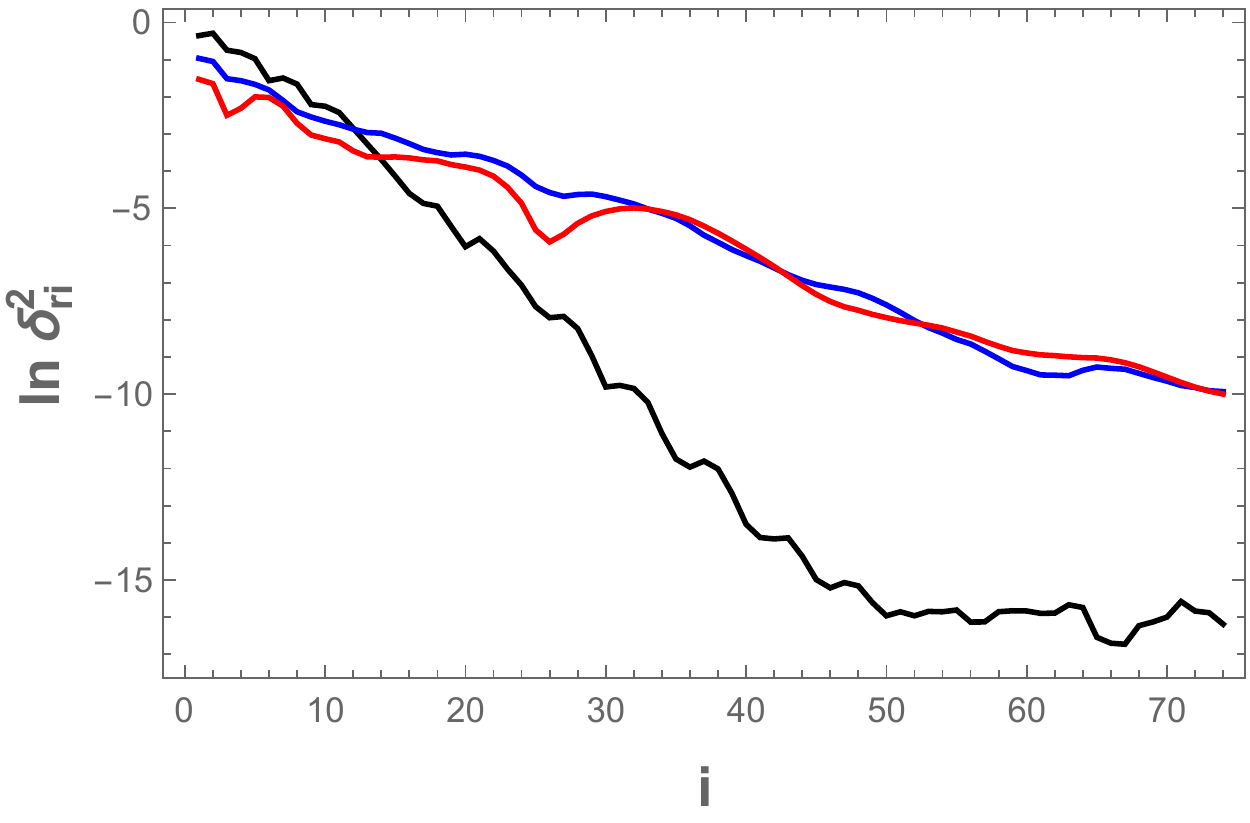,width=7.5cm}
& 
\psfig{figure=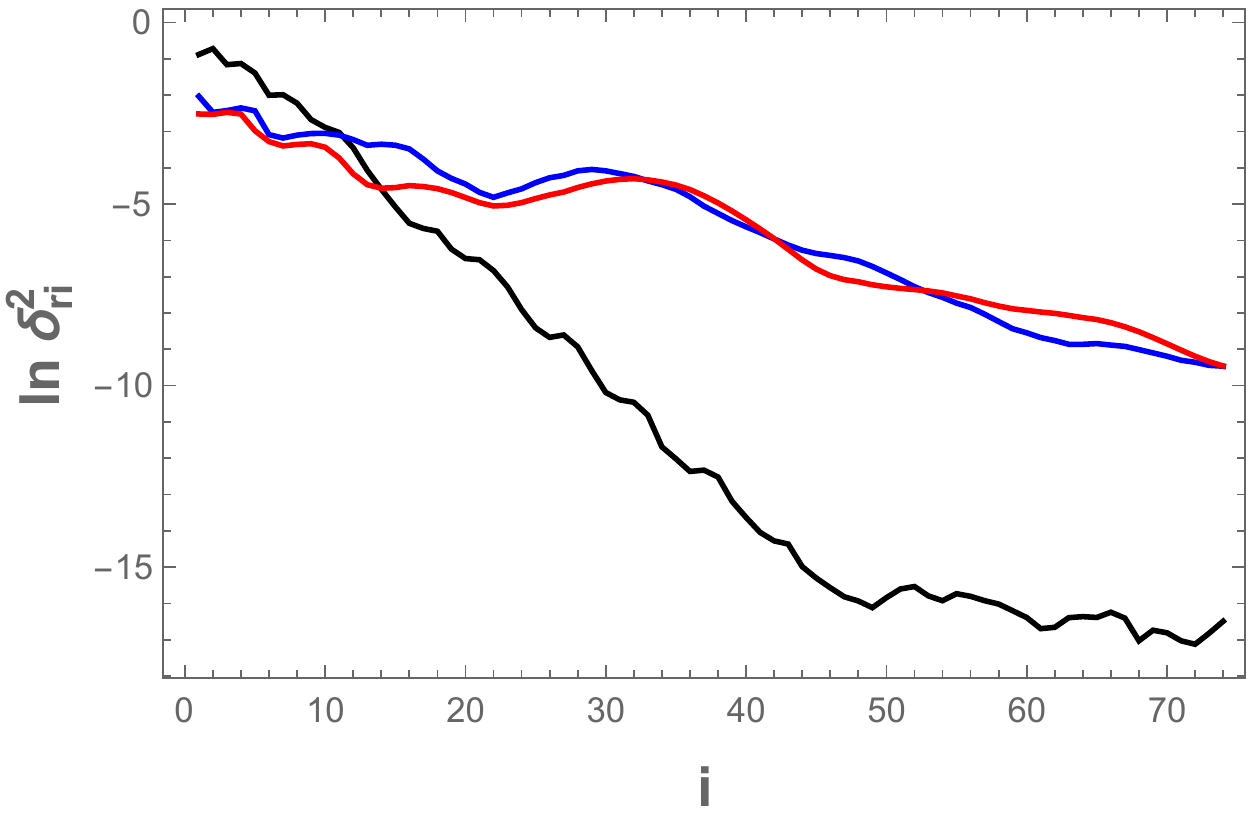,width=7.5cm} 
\\
{\bf a}  
& {\bf b} 
\end{tabular}
\caption{The dynamics of the shifted difference norm 
$\delta_{ri}^2$  (shown in logarithmic scale) in the $r$-th RN of the chain for
$r=1$ (black), $r=2$ (blue), $r=3$ (red) and $m=75$.
({\bf a}) The MW algorithm with renewal of the state vectors
$\bm s_{r0}^j$ before start of the prediction round.
({\bf b}) The MW algorithm with inheritance of the state vectors
$\bm s_{r0}^j = \bm s_{rm}^{j-1}$ between the rounds.
}
\label{delta}
\end{center}
\end{figure} 

First we consider the shifted difference dynamics for the MW algorithm
assuming that at the beginning of each prediction round the 
initial RN state is refreshed by setting it to zero $\bm s_{r0}^j = \bm 0$. Fig. \ref{delta}a
shows that the state vectors $\bm s_{ri}^j$ for $j=2$ 
and $r=1,2,3$ decay exponentially. 
The largest decay rate is observed in
the first RN ($r=1$) and the norm $\delta_{1i}^2$ levels off
for large $i > 50$. Nevertheless, for two other RNs the decay
is also significant and the norm $\delta_{rm}^2$ is negligibly small
compared to the characteristic norm of the state vector itself
$\delta_{r,m-1}^2 \ll s_{rm}^2, s_{r,m-1}^3 \sim 1$.

Similar behavior we observed in case when at the prediction round $j$ the 
initial state $\bm s_{r0}^j$ of the $r$-th RN is inherited from the 
previous round, namely, $\bm s_{r0}^j = \bm s_{rm}^{j-1}$. The results are shown in 
Fig. \ref{delta}b.
Assuming that the results presented in Fig. \ref{delta} remain valid for all
prediction rounds we start development of the ML algorithm generalization
for RN chains.

\begin{figure}[h!]
\begin{center}
\psfig{figure=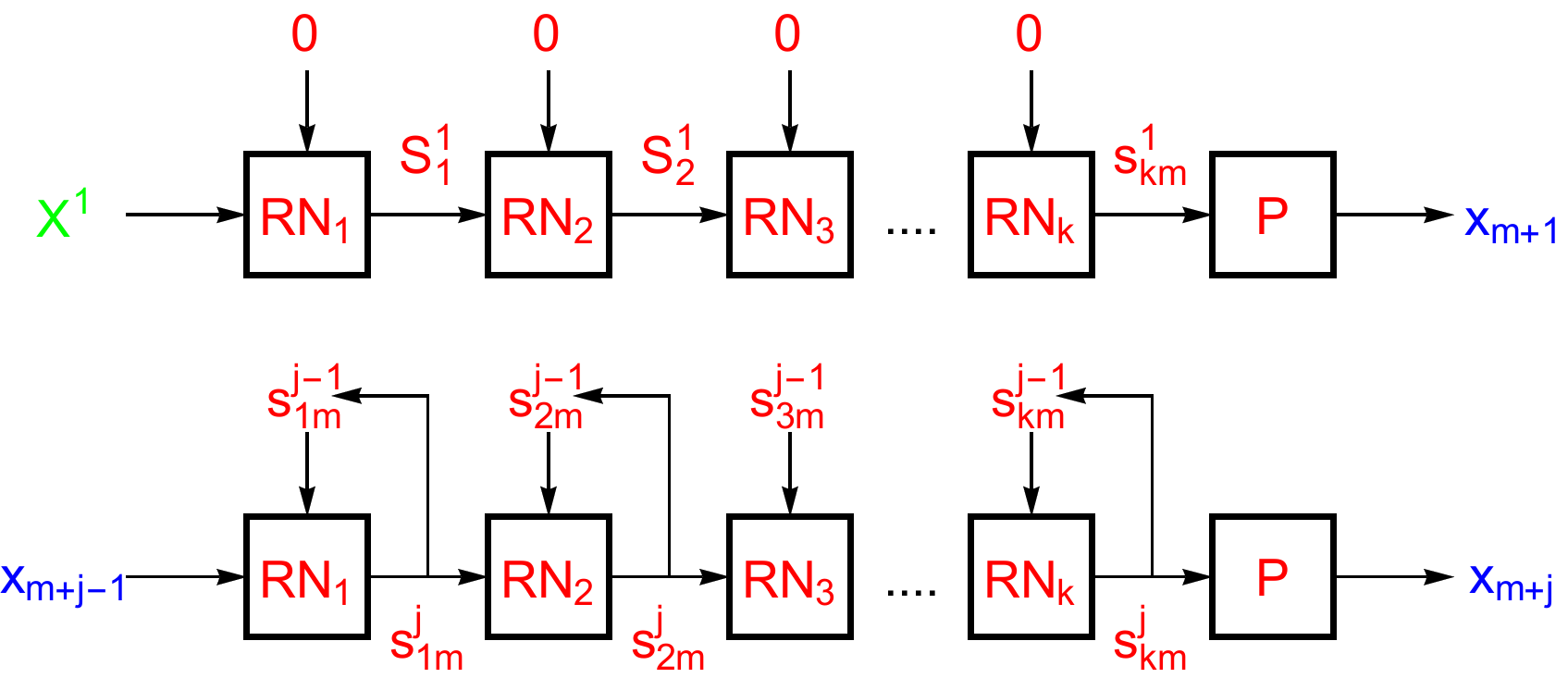,width=14cm}
\caption{The predictive ML algorithm in the chain of RN${}_r$
followed by the predictor P for the first  (top panel) and 
subsequent (bottom panel) predictive rounds.
The first round is similar to the one in the MW algorithm while
all subsequent rounds use the transformations (\ref{state_last_r}, \ref{state_last_1}).
The new states $\bm s_{rm}^{j}$ are used as the corresponding RN states in the next round.
}
\label{ML}
\end{center}
\end{figure}

Note that for the $r$-th RN the dynamics of its state governed by
(\ref{F_r}) which gives for $i=m$ at the $j$-th prediction round ($j>1$)
\be
\bm s_{r,m}^j = \bm F_r(\bm s_{r-1,m}^j, \bm s_{r,m-1}^j) = 
\bm F_r(\bm s_{r-1,m}^j, \bm s_{r,m}^{j-1}).
\label{state_last_r}
\ee
This transformation relates the final $m$-th states
$\bm s_{r,m}^j,\ \bm s_{r,m}^{j-1}$ of 
the $r$-th RN at two adjacent predictive rounds and the 
final $m$-th state $\bm s_{r-1,m}^j$ of the preceding $(r-1)$-th RN
at the current $j$-th round.
The relation for first ($r=1$) RN reads
\be
\bm s_{1,m}^j = \bm F_1(\bm s_{0,m}^j, \bm s_{1,m-1}^j) = 
\bm F_1(\bm x_{m+j-1}, \bm s_{1,m}^{j-1}) = 
\bm F_1(\bm p(\bm s_{k,m}^{j-1}), \bm s_{1,m}^{j-1}).
\label{state_last_1}
\ee
These two relations completely determine the dynamics of the 
final states $\bm s_{r,m}^j$ for $1 \le r \le k$ and $2 \le j \le p$.
The computation of the states $\bm s_{r,m}^1$ and the first
predicted value $\bm x_{m+1}$ is performed 
using once the standard MW algorithm. The corresponding 
schematics is depicted in Fig. \ref{ML}.
The second prediction round starts with feeding the first
predicted value $\bm x_{m+1} = \bm s_{0,m}^2$ into the first RN
being in the state $\bm s_{1m}^1$ inherited from the first round.
The result of the transformation (\ref{state_last_1}) is a new state
$\bm s_{1m}^2$ of the first RN which in its turn represents the input fed into the
second RN. At this moment the state of the second RN is $\bm s_{2m}^1$
and the transformation (\ref{state_last_r}) produces a new state
$\bm s_{2m}^2$. This procedure is repeated for all remaining RNs in the 
chain until the last $k$-th RN is reached. The result of (\ref{state_last_r})
application is $\bm s_{km}^2$ that is used by the predictor P to generate
a new prediction $\bm x_{m+2}$ which serves as an input to the next round of 
prediction (bottom panel in  Fig. \ref{ML}). It is reasonable to assume that an
implementation of the memoryless algorithm is innate in brain environment.

The total number of transformations $N_{ML}^k = mk+1+(k+1)(p-1) = k(m+p-1)+p$ compared to 
$N_{MW}^k = (mk+1)p$ leads to the speed gain 
\be
\gamma_k=\frac{N_{MW}^k}{N_{ML}^k} = 
\frac{(mk+1)p}{k(m+p-1)+p}.
\label{gammak}
\ee

To address the problem of prediction quality of ML approach compared to 
MW algorithm we perform computation of specified number $p$ of the predicted values
for the different length $m$ of the input sequence using both algorithms.
 
\begin{figure}[h!]
\begin{center}
\begin{tabular}{cc}
\psfig{figure=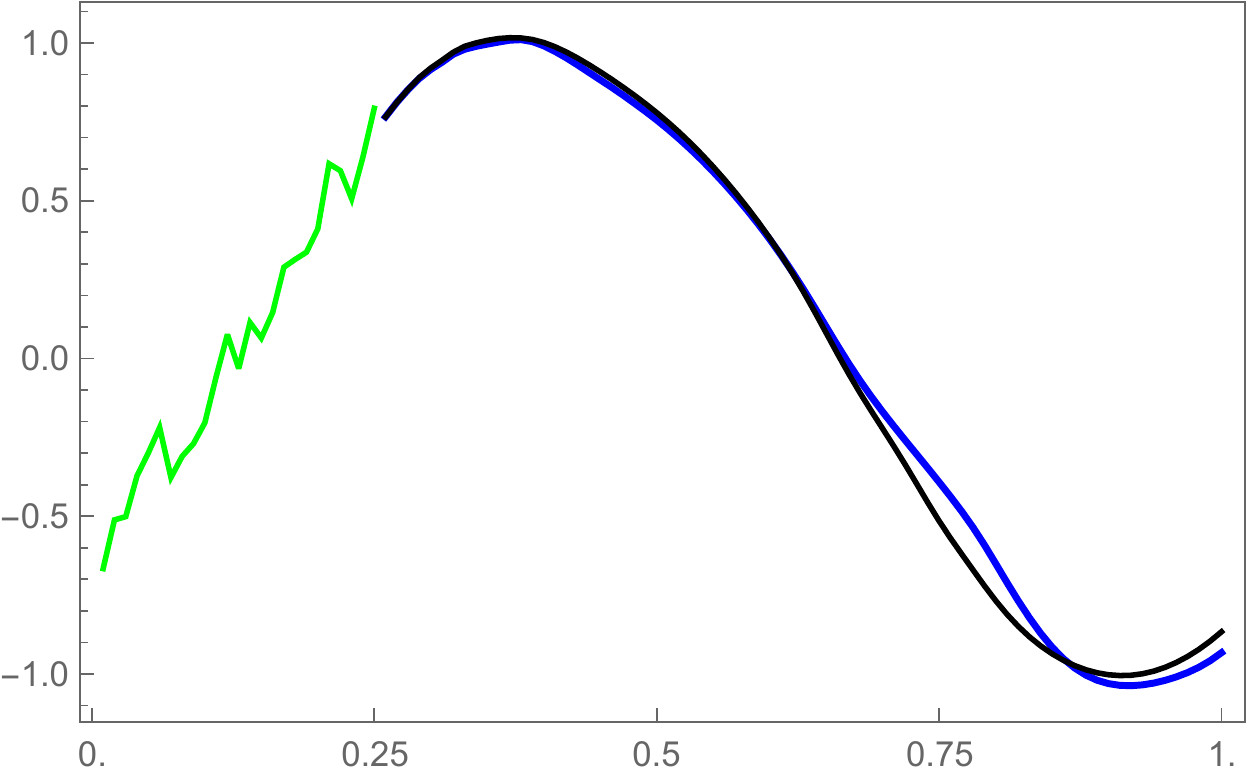,height=4.5cm}  &
\psfig{figure=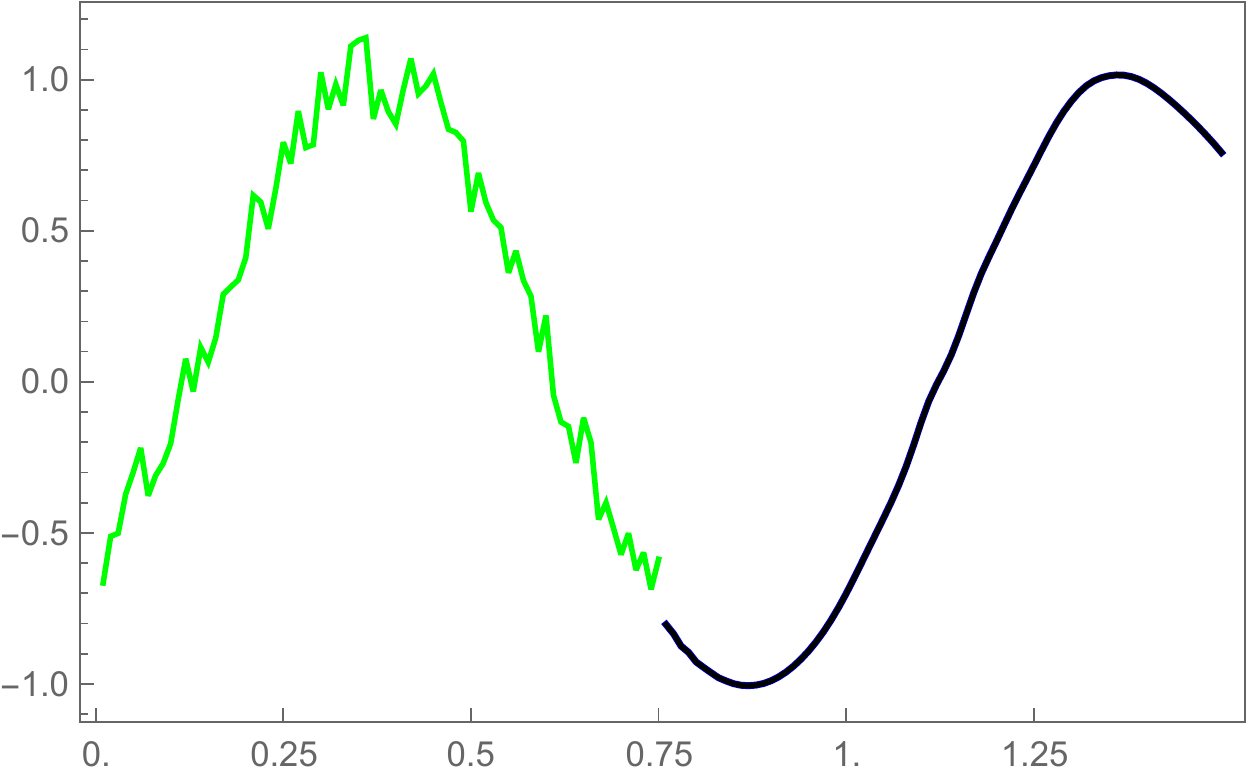,height=4.5cm}
\\
{\bf a} & {\bf b}
\end{tabular}
\caption{Comparison of the predictions  of $p=75$ elements by the moving window (blue) 
and the memoryless (black) algorithms for the sine wave with noise amplitude $a=0.15$ input
sequence (green) of length ({\bf a}) $m=25$ and ({\bf b}) $m=75$;
in ({\bf b})  the curves coincide.}
\label{MWMLcomparison}
\end{center}
\end{figure}

The results presented in
Fig. \ref{MWMLcomparison} show that the increase in $m$ makes ML prediction
to coincide with that of by MW algorithm. It can be explained by inspection of the 
shifted difference behavior (see Fig. \ref{delta}) -- for large $m$ the 
value of $\bm \delta_{rm}^j$ is negligibly small that allows the approximate
transformations in (\ref{state_last_r},\ref{state_last_1}) more precisely reproduce the 
prediction dynamics determined by the MW approach.

Thus it is reasonable to consider cases where the number $p$  of predicted points 
is usually less or approximately equal to the length $m$
of the initial input sequence. Then the speed gain estimates is in the range from
$\gamma_k \sim p$ for $p \ll m$ to $\gamma_k \sim mk/(2k+1)$ for $p \sim m$.
For the less frequent case $p \gg m$ the speed gain is 
$\gamma_k \sim mk/(k+1)$.

\begin{figure}[h!]
\begin{center}
\begin{tabular}{cc}
\psfig{figure=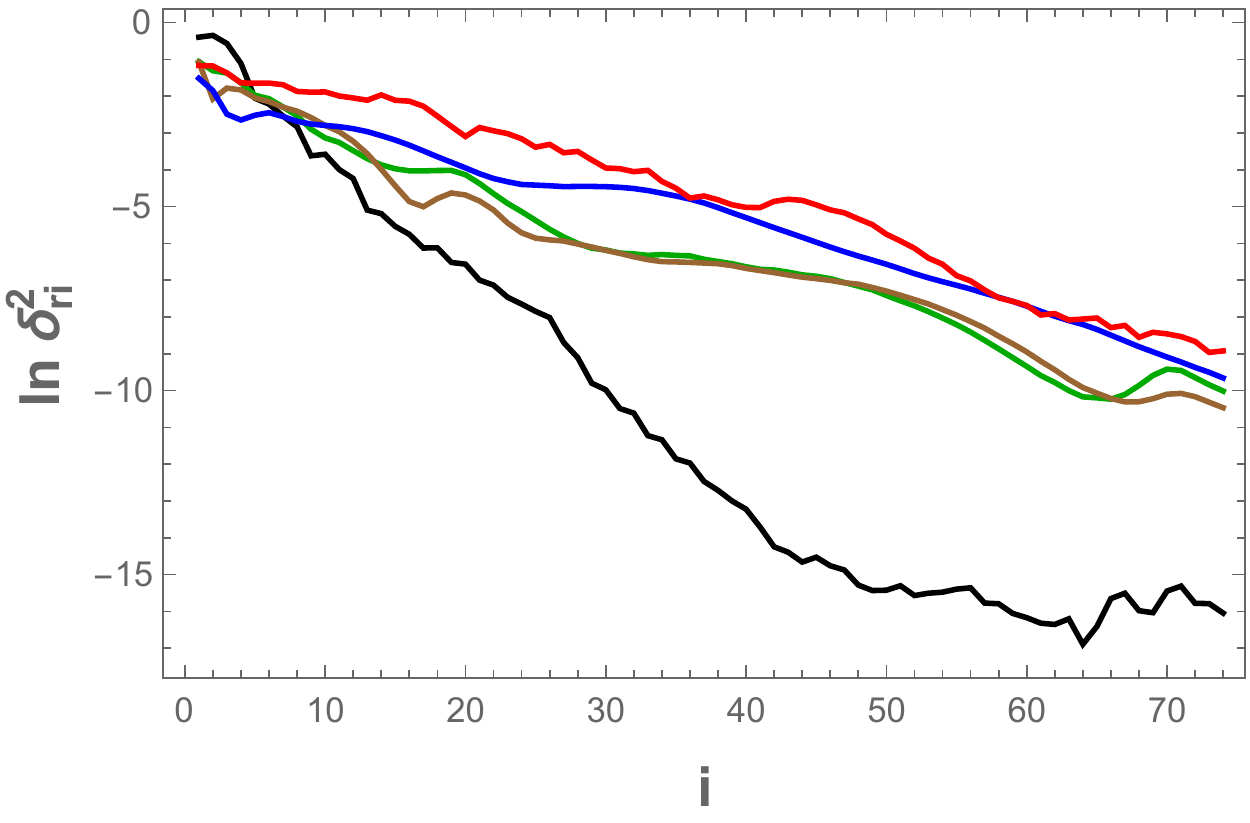,width=7.5cm}
& 
\psfig{figure=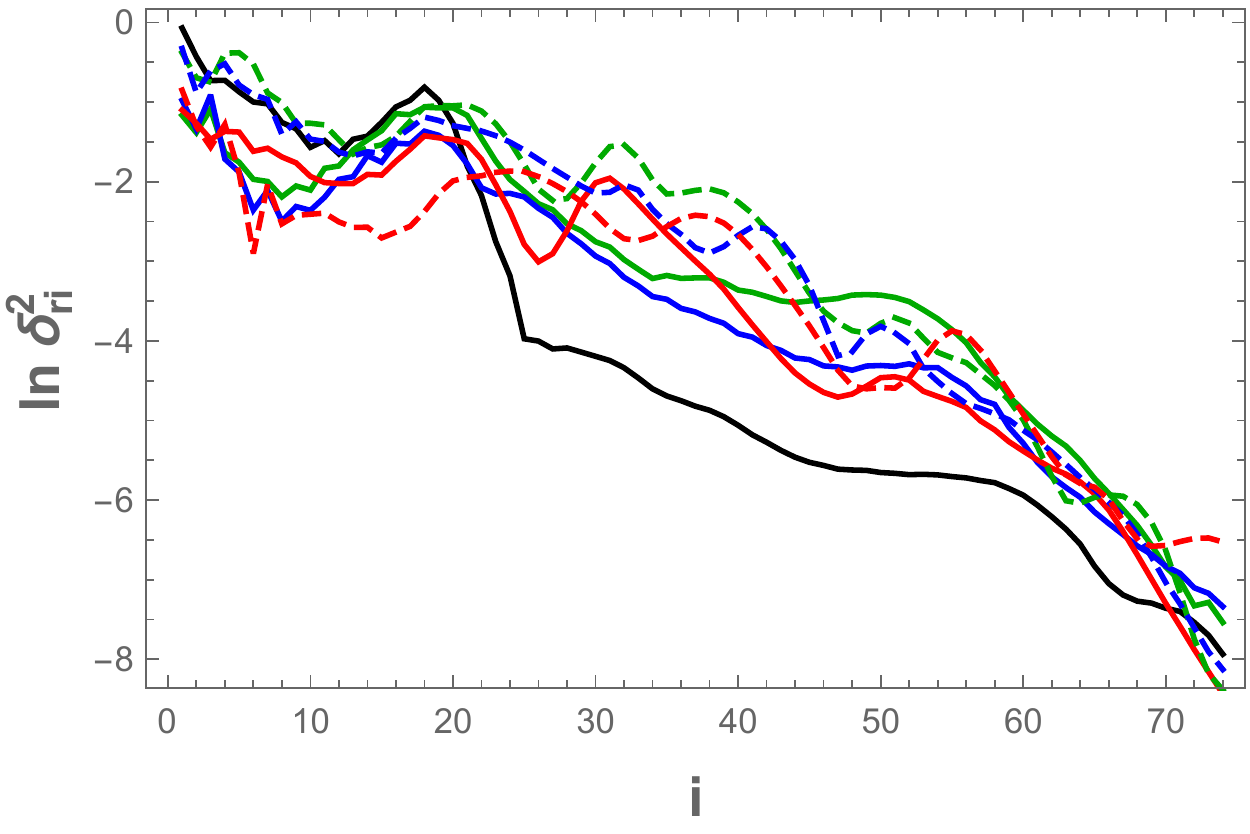,width=7.5cm} 
\\
{\bf a}  
& {\bf b} 
\end{tabular}
\caption{The dynamics of the shifted difference norm
$\delta_{ri}^2$  in the $r$-th RN of the chain with $m=75$
and ({\bf a}) $k=5$ and ({\bf b}) $k=7$.
The black curve corresponds to $r=1$.
}
\label{delta5_7}
\end{center}
\end{figure}

It is instructive to consider the RN chains with increased number of recurrent modules.
We tested the networks with $k=5$ and $k=7$ and found that the 
dynamics of the shifted difference norm demonstrates a lower rate of exponential decay
(Fig. \ref{delta5_7}) for the larger number of RNs.
The same time ML algorithm for $k=5,7$
continues to produce predicted trajectories of the same quality as MW approach does 
(similar to those shown in Fig. \ref{MWMLcomparison}).

\section{ML algorithm for well trained network}
\label{well_trained}
The numerical experiments discussed in the previous section 
show that the high quality of prediction makes possible the existence of
the ML algorithm even when the the shifted difference decay is not strongly 
pronounced in some of (or all) RN modules. 
To simplify the presentation of the results we use in this section a 
different notation for the main transformation (\ref{F_r}), namely
\be
\bm s_{r,i} = \bm F_r(\bm s_{r-1,i}, \bm s_{r,i-1}) \equiv
\bm F_r[\bm s_{r,i-1}](\bm s_{r-1,i}).
\label{F_r_new}
\ee
The relations (\ref{state_last_r}) and (\ref{state_last_1}) take form
\be
\bm s_{1,m}^{j+1} =
\bm F_1[\bm s_{1,m}^{j}](\bm p(\bm s_{k,m}^{j})).
\quad\quad
\bm s_{r,m}^{j+1} =
\bm F_r[\bm s_{r,m}^{j}](\bm s_{r-1,m}^{j+1}).
\label{state_transform}
\ee

First consider a network with a single RN module and a predictor
discussed in Section \ref{single}.
Note that the condition (\ref{last_state_diff}) implies 
(\ref{ML_transform}) for the MW algorithm independent 
of the prediction quality. In other words, the network 
can be trained badly (for example, due to a small number of neurons $n$)
but if  (\ref{last_state_diff}) holds the network prediction
will be governed by  (\ref{ML_transform}).

On the other hand, in case when the network is trained well
(i.e., the norm of the difference $|\bar{\bm x}_{m+1}-\bm x_{m+1}|$ between the predicted value and
the ground truth is negligible compared to the characteristic range $A \sim 1$ of the 
sequence $\bm X$)
the application of EW algorithm leads to this transform naturally.
Consider a recurrent network of general type which state dynamics is governed by 
(\ref{F}) and use it for prediction of the value $\bar{\bm x}_{m+j+1}$
based on the input sequence $\bm X_{m+j}=\{\bm x_i\},\ 1\le i\le m+j,$ consisting of $m+j$ elements
$\bm x_i$. The prediction $\bar{\bm x}_{m+j+1} = \bm p(\bm s_{m+j})$ 
 is performed by applying the transformation $\bm p$ to the 
final RN state $\bm s_{m+j}$.
Use the same RN to predict the next point $\bar{\bm x}_{m+j+2}$ using the expanded
input sequence of $m+j+1$ elements  $\bm X_{m+j+1}=\{\bm x_i\},\ 1\le i\le m+j+1$ and
find the final network state
$\bm s_{m+j+1} = \bm F[\bm s_{m}](\bm x_{m+j+1})$.
Assuming that the prediction quality of the RN is good enough
replace in the above relation $\bm x_{m+j+1}$ by $\bar{\bm x}_{m+j+1}$ 
as these two values are close to each other and obtain
$$
\bm s_{m+j+1} \approx \bm F[\bm s_{m+j}](\bar{\bm x}_{m+j+1}) =
\bm F[\bm s_{m+j}](\bm p(\bm s_{m+j})) \equiv 
F(\bm p(\bm s_{m+j}), \bm s_{m+j}),
$$
which is equivalent  to (\ref{ML_transform}). Note that in this case
the condition $\epsilon^j \ll 1$ in (\ref{last_state_diff}) is not
required. 

Turn to the predictive network made of $k$ RNs and a predictor
considered in details in Section \ref{multi}. For the 
input sequence $\bm X_{m+j}=\{\bm x_i\},\ 1\le i\le m+j,$ we predict 
$\bar{\bm x}_{m+j+1} = \bm p(\bm s_{k,m+j})$
which by assumption is close to last 
element $\bm x_{m+j+1}$ of the expanded
input sequence of $m+j+1$ elements  $\bm X_m=\{\bm x_i\},\ 1\le i\le m+j+1$.
Use $\bm X_{m+j+1}$ and find for $r=1$
\be
\bm s_{1,m+j+1} =  
\bm F_1[\bm s_{1,m+j}](\bm x_{m+j+1}).
\label{state_last_1EW}
\ee
Use (\ref{F_r_new}) to write for $r = 2$:
\be
\bm s_{2,m+j+1} = \bm F_2[\bm s_{2,m+j}](\bm s_{1,m+j+1}) = 
\bm F_2[\bm s_{2,m+j}]\circ\bm F_1[\bm s_{1,m+j}](\bm x_{m+j+1}),
\label{state_last_1EW}
\ee
where $\circ$ denotes function composition: $R\circ P(x) \equiv R(P(x))$.
Continue for $r > 2$ to obtain
\bea
\bm s_{3,m+j+1} &=& \bm F_3[\bm s_{3,m+j}](\bm s_{2,m+j+1}) = 
\bm F_3[\bm s_{3,m+j}]\circ\bm F_2[\bm s_{2,m+j}]\circ\bm F_1[\bm s_{1,m+j}](\bm x_{m+j+1}),
\nonumber \\
& \ldots &
\nonumber \\
\bm s_{k,m+j+1} &=& \bm F_k[\bm s_{k,m+j}](\bm s_{k-1,m+j+1}) = 
\bm F_k[\bm s_{k,m+j}]\circ \ldots \circ\bm F_2[\bm s_{2,m+j}]\circ\bm F_1[\bm s_{1,m+j}](\bm x_{m+j+1}).
\label{EW_chain}
\eea
Now in the last equation in (\ref{EW_chain}) replace $\bm x_{m+j+1}$ by $\bar{\bm x}_{m+j+1} = \bm p(\bm s_{k,m+j})$
and find
\be
\bm s_{k,m+j+1} = 
\bm F_k[\bm s_{k,m+j}]\circ \ldots \circ\bm F_2[\bm s_{2,m+j}]\circ\bm F_1[\bm s_{1,m+j}]( \bm p(\bm s_{k,m+j})).
\label{state_last_kEW}
\ee
It is easy to see that (\ref{state_last_kEW}) is equivalent to the relations
(\ref{state_transform}), this equivalence implies that the well trained network 
represented by a chain of RNs can be used for prediction by
memoryless algorithm.

\section{Discussion}
\label{discussion}

In this manuscript we consider a general type predictive network 
that transforms an input sequence made of elements with similar structure to predicts a
single element of the same type. 
The input sequence flows first into an encoder then to a chain of the recurrent networks (RNs)
to a predictor and finally to a decoder.
We focus on the network central part (the RN chain and the predictor) responsible for
the actual prediction procedure. 
When one needs to generate several predicted values there exists the
traditional moving window (MW) predictive algorithm well fitted
for usage in the artificial neural networks. However its
implementation in the natural neural network existing in brain can 
fail as some requirement are very difficult to satisfy and the prediction
loses both robustness and reliability. 
In \cite{Rub2020a} for the simplest predictive network 
we proposed another faster and more stable
algorithm that that does not require memorization of the input sequence
(as a whole or a part of it) and we call it {\it memoryless} (ML) algorithm.
We show that the ML approach can be generalized to the case of the 
RN chain (compared of a single RN in \cite{Rub2020a}).

The essence of the approach is that it requires the 
input sequence only once to predict the first value and to update the 
state of each RN in the chain. 
Then the original input can be forgotten and consecutive prediction rounds
are based on the own network dynamics. The input sequence is replaced by the 
latest predicted value which is transformed in the network to 
update the RN states and generate the next prediction.

We show that ML algorithm produces results that are very close to those generated by 
MW approach. The increase of the input sequence length makes 
ML results to coincide with MW ones but gives the speed up proportional to 
the number of predicted points. We also discuss a possibility of the ML algorithm
implementation in brain networks and show that it is quite suitable 
to this purpose.

Based on the MW algorithm analysis we establish a set of conditions
for ML application irrespective of the predictve quality of the network that
demand the shifted difference norms $\delta_{r,i}^j$ to decay exponentially. 
It appears that when the network is well trained the conditions on the difference norms
can be dropped. It is quite important result as in the natural neural systems 
the analysis of such behavior is very difficult if not just impossible.
On the other hand the evolution of natural neuron networks should select
only those networks that have highest predictive ability which in its turn
successfully implement the robust and fast memoryless algorithm.
This observation makes the proposed ML algorithm a critical element
of a self-consistent scheme based on neural networks of high performance -- 
it allows to achieve fast and extremely robust predictions with these networks.

The number of transformation layers in a chain is an important factor
determining both complexity and prognostic ability of a predictive network.
For example, the GPT (Transformer) network that uses more complex (compared to (\ref{F})) attention
signal transformation \cite{Vaswani2017} demonstrates
an increasing ability to solve successfully complex language-based tasks 
when the number of elements in the chain grows from 12 to 96
with parallel increase of total number of neurons in individual attention block
from 768 to 12288 \cite{Brown2020}. 
In brains with small total number of neurons the number $k$ of RNs in the chain might be limited
and it can be ascribed to a restriction of the total number of neurons
assigned for a specific networks. It is instructive to know how 
increase in the RN number affects possibility of successful implementation 
of ML algorithm especially in case when the condition
of the exponential decay of the shifted difference norm $\delta_{r,i}^j$ fails.
The answer to this question requires to train RN chains with 
increasing $k$ and test whether the shifted difference exponential decay
is preserved. We tested the chains of $k=3,5,7$ RNs and found that
the decay rate decreases with increase of $k$. It is reasonable to 
extend such analysis to much larger number (several dozens) of RN modules in the chain.

\section*{Acknowledgements}
The author wishes to thank Yuri Shvachko, Kausik Si  and Mitya Chkolvsky for fruitful discussions.

\end{document}